\documentclass[12pt,eqno]{article}
\usepackage{amsmath,amsthm,amscd,amssymb}
\usepackage{latexsym}
\newtheorem{remark}{Remark}[section]
\newcommand{\bremark}{\begin{remark} \em}
\newcommand{\eremark}{\end{remark} }

\begin{document}
\parindent 12pt
\renewcommand{\theequation}{\thesection.\arabic{equation}}
\renewcommand{\baselinestretch}{1.15}
\renewcommand{\arraystretch}{1.1}
\def\disp{\displaystyle}
\title{\large Iterative method for Kirchhoff-Carrier type equations and its applications\\
\footnotetext{E-mail addresses:qiuyidai@aliyun.com.\\
This work is supported by National Natural Science Foundation of China (Grant:No.11671128).}}
\author{{\small Qiuyi Dai}\\
{\small\emph{ College of Mathematics and Statistics,Hunan Normal University,}}\\
\small\emph{{Changsha Hunan 410081,PR China}}\\
\date{}}
\maketitle
\abstract
{\small Let $A(s, t)$ be a continuous function with a positive lower bound $m$, and $\Omega$ be a bounded domain in $R^N$. In this short note, we use a so-called invariant set for finding nonnegative solutions of the following Kirchhoff-Carrier type equations
\begin{equation*}
\left\{\begin{array}{ll}
-A(\|u\|_p, \|\nabla u\|_2)\Delta u=g(x, u) & x\in\Omega,\\
u=0 & x\in\partial\Omega.
\end{array}
\right.
\end{equation*}
The main advantage of our method is that we almost need no restrictions on $A(s, t)$ except for continuous and a positive lower bound. This removes away the monotonicity assumption of $A(s, t)$ used in most papers based on sub-supersolution method. As applications of the abstract result obtained by our method, some concrete examples are also studied in Section 2 of this paper.}

\section*{1. Invariant set method}

\setcounter{section}{1}

\setcounter{equation}{0}

\noindent

Let $A:\ R_+\times R_+\to R$ be a continuous function with $A(s, t)\geq m$ for some positive constant $m$. If necessary, $A(s, t)$ can be replaced by $A(s, t)/m$. Therefore, no loss of generality, we may assume $m\geq1$. For $a, b\geq 0$, we set
$$M_{a,b}=\max\{A(s, t):\ (s, t)\in [0, a]\times[0, b]\}.$$

Let $\Omega$ be a bounded domain in $R^N$. Denote by $L^p(\Omega)$ the standard Lebesgue's space with norm $\|\cdot\|_p$. Let $H_0^1(\Omega)$ be the standard Sobolev space with norm $\|\nabla\cdot\|_2$, and $2^*$ be the critical exponent for the Sobolev embedding (that is $2^*=+\infty$ if $N\leq 2$, and $2^*=2N/(N-2)$ if $N\geq 3$). Let $\Delta$ be the Laplace operator, and $g(x,s)$ be a H\"older continuous fuction defined on $\overline{\Omega}\times R$. We consider the following problem of generalized Kirchhoff-Carrier type equations 
\begin{equation}\label{eq3}
\left\{
\begin{array}{ll}
-A(\|u\|_p, \|\nabla u\|_2)\Delta u=g(x, u) & x\in\Omega,\\
u=0& x\in\partial\Omega.
\end{array}
\right.
\end{equation}
Assume that $g(x, s)$ satisfies.
\vskip 0.1in
{\bf (G1):}\ There are two functions $0\leq\varphi(x)\leq\psi(x)$ such that
\begin{equation}\label{eq1}
\left\{\begin{array}{ll}
-\Delta\varphi\leq g(x,\varphi) & x\in\Omega,\\
-\Delta\psi\geq g(x,\psi) & x\in\Omega,\\
\psi\geq\varphi=0 & x\in\partial\Omega.
\end{array}
\right.
\end{equation}
\vskip 0.1in
{\bf (G2):}\ There exists a real number $\alpha$ such that for any $0<\beta\leq 1$ and $\omega\in [\beta\varphi, \psi]$, there hold
$$-\beta^\alpha\Delta\varphi\leq g(x,\omega)\leq-\Delta\psi.$$

Our model problem (\ref{eq3}) covers the classical Kirchhoff problem which corresponds to the case $A(\|u\|_p, \|\nabla u\|_2)=1+d\|\nabla u\|_2^2$ (see \cite{Kir}), and the classical Carrier problem which corresponds to the the case $A(\|u\|_p, \|\nabla u\|_2)=1+d\|u\|_2^2$. The Kirchhoff problem has variational structure, and hence can be studied by powerful variational method. Therefore, it has attracted many attentions in recent years (see for example \cite{HZ, LLS, L, AZ, MPR, Na1, Na2, PZ, YRA, ZSN, ZS}). The Carrier problem has no variational structure, and less literatures on carrier problem can be found (see however \cite{CDCS, CL, CR, C, ALP}). In this short note, we propose an invariant set method for finding solutions of problem (\ref{eq3}). Iterative procedure based on comparison principle of Kirchhoff type operator itself has been used to find solution to Kirchhoff equations by many authors (see for example \cite{AC1, AC2, HD, DM, ACS, FS, Ng}). However, comparison principle may cease to validate for general Kirchhoff type operator except for that possessing some monotonicity property. Instead, our invariant set method based only on the comparison principle of Laplace operator. Therefore, we only need to put mini restrictions on Kirchhoff type operator. To state our main result, we let $M=\max\{A(s, t):\ (s, t)\in [0, \|\psi\|_p]\times [0, \|\nabla\psi\|_2]\}$ and
\begin{equation}\label{req 1}
r_M=\left\{\begin{array}{ll}
1&\mbox{if}\ M=1,\\
(\frac{1}{M})^{\sum\limits_{k=0}^\infty\alpha^k}&\mbox{if}\ M>1\ \mbox{and}\ \sum\limits_{k=0}^\infty\alpha^k<+\infty,\\
0&\mbox{if}\ M>1\ \mbox{and}\ \sum\limits_{k=0}^\infty\alpha^k=+\infty.
\end{array}
\right.
\end{equation}
It is easy to see that $\frac{1}{M}r_M^\alpha=r_M$. By making use of invariant set and Schaulder fixed point theorem, we can prove the following result.
\vskip 0.1in
{\em {\bf Theorem 1.1}\  Assume that $1<p\leq 2^*$, and {\bf (G1)-(G2)} hold. Then problem (\ref{eq3}) has at least one nonnegative solution $u(x)$ with property $r_M\varphi(x)\leq u(x)\leq\psi(x)$ for any $x\in\Omega$. }
 \vskip 0.1in
 {\bf Proof:}\ Let $C(\overline{\Omega})$ be the Banach space consisting of all continuous functions defined on $\overline{\Omega}$ and endowing with norm 
 $\|u\|_C=\max\{|u(x)|:\ x\in\overline{\Omega}\}$ for any $u\in C(\overline{\Omega})$. Setting $X=H_0^1(\Omega)\cap C(\overline{\Omega})$ and endowing $X$ with norm $\|u\|_X=\|\nabla u\|_2+\|u\|_C$ for any $u\in X$, it is easy to check that $(X, \|\cdot\|_X)$ is a Banach space. Let $G(x, y)\geq 0$ denote the Green's function on $\Omega$. We define a map $T$ on $X$ by 
 $$T(v)(x)=\frac{1}{A(\|v\|_p, \|\nabla v\|_2)}\int_\Omega G(x, y)g(y, v(y))dy\quad\quad\forall v\in X.$$
 Obviously, $T$ is continuous on $X$. Moreover, the following two claims hold for map $T$.
\vskip 0.05in 
 {\bf Claim (i):}\ If we set $[r_M\varphi, \psi]=\{v(x)\in X:\ r_M\varphi(x)\leq v(x)\leq\psi(x)\ \forall x\in\Omega\}$, then we can claim that $[r_M\varphi, \psi]$ is an so-called invariant set of map $T$. That is
 $$T([r_M\varphi, \psi])\subset [r_M\varphi, \psi].$$
 
 In fact, for any $w\in [r_M\varphi, \psi]$, it follows from {\bf (G2)} that
 $$
 \begin{array}{ll}
 T(w)(x)&=\frac{1}{A(\|w\|_p, \|\nabla w\|_2)}\int_\Omega G(x, y)g(y, w(y))dy\\
 &\geq\frac{1}{A(\|v\|_p, \|\nabla v\|_2)}r_M^\alpha\int_\Omega G(x, y)(-\Delta\varphi(y))dy\\
 &=\frac{1}{A(\|v\|_p, \|\nabla v\|_2)}r_M^\alpha\varphi(x)\\
 &\geq\frac{1}{M}r_M^\alpha\varphi(x)=r_M\varphi(x). 
 \end{array}
 $$
and 
$$
\begin{array}{ll}
T(w)(x)&=\frac{1}{A(\|w\|_p, \|\nabla w\|_2)}\int_\Omega G(x, y)g(y, w(y))dy\\
 &\leq\frac{1}{A(\|v\|_p, \|\nabla v\|_2)}\int_\Omega G(x, y)(-\Delta\psi(y))dy\\
 &=\frac{1}{A(\|v\|_p, \|\nabla v\|_2)}\psi(x)\\
 &\leq\psi(x). 
 \end{array}
 $$
 Therefore, $T(w)(x)\in [r_M\varphi, \psi]$. This concludes Claim (i).
 \vskip 0.05in
 {\bf Claim (ii):}\ Let $G=\max\{g(x, s):\ (x, s)\in\overline{\Omega}\times[0, \|\psi\|_C]\}$. Then, for any $w\in [r_M\varphi, \psi]$ and $\tau\in (0, 1)$, $T(w)(x)$ is in $C^{1,\tau}(\Omega)$ with $C^{1,\tau}(\Omega)$ being the standard H\"{o}lder space. Moreover, there exist a positive constant $C_1$ depending only on $N$ and $\Omega$ such that
 $$\|T(w)\|_{C^{1,\tau}(\Omega)}\leq C_1G\quad\quad\forall w\in[r_M\varphi, \psi].$$ 
  
  In fact, by the definition of $G$ and the assumption of $A(s, t)$, we have
  $$|\frac{g(x, w(x))}{A(\|w\|_p, \|\nabla w\|_2)}|\leq G\quad\quad\forall w\in [r_M\varphi, \psi].$$
  Therefore, the conclusions of Claim (ii) follow from properties of Green's function which can be found in the textbook \cite{GT} of Gilbarg and Trudinger.
  
  Let $D=\{v\in X:\ v(x)\in[r_M\varphi, \psi]\ \mbox{and}\ \|v\|_{C^{1,\tau}(\Omega)}\leq C_1G\}$. Then, it is easy to see that $D$ is convex and closed. Moreover, by Ascoli-Arzela's theorem, we know that $D$ is sequential compact in $X$. Therefore $D$ is a compact set of $X$. By Clam (i) and (ii), we know that $T$ maps $D$ into $D$. It follows from the continuity of $T$ and the Schaulder fixed point theorem that there exists a function $u(x)\in D$ such that
  $$u(x)=\frac{1}{A(\|u\|_p, \|\nabla u\|_2)}\int_\Omega G(x, y)g(y, u(y))dy.$$
  By the regularity theory of elliptic equations, we can see that $u(x)$ is a classical solution to problem (\ref{eq3}) with property $r_M\varphi(x)\leq u(x)\leq\phi(x)$ for any $x\in\Omega$. This completes the proof of Theorem 1.1.

\section*{2. Applications}

\setcounter{section}{2}

\setcounter{equation}{0} 

\noindent
 
As applications of Theorem 1.1, we give some concrete examples in this section. The first example is a Drichlet problem of inhomogeneous Carrier's equation.
\vskip 0.1in
{\bf Example 1:}\ Let $f(x)\in C^1(\Omega)\backslash\{0\}$ be a function such that the following problem has a solution
\begin{equation}\label{eq21}
\left\{
\begin{array}{ll}
-\Delta\phi=f(x) & x\in\Omega,\\
\phi\geq 0 & x\in\Omega,\\
\phi=0 & x\in\partial\Omega.
\end{array}
\right.
\end{equation}
Obviously, $f(x)$ may change sign. For $d, p>0$, we consider the following problem with positive parameter $\lambda$.
\begin{equation}\label{eq22}
\left\{
\begin{array}{ll}
-(1+d\|u\|_2^2)\Delta u=u^p+\lambda f(x) & x\in\Omega,\\
u>0 & x\in\Omega,\\
u=0 & x\in\partial\Omega.
\end{array}
\right.
\end{equation}
It is worth pointing out that Carrier's operator $-(1+d\|\cdot\|_2^2)\Delta$ is not a variational operator, and it also lack of comparison principle. Therefore, powerful tools of variational and sub-supersolution method can not be used directly to study problem (\ref{eq22}). Here, we use Theorem 1.1 to conclude the following result.
\vskip 0.1in
{\em {\bf Theorem 2.1}\ There exists a positive number $\lambda_f$ such that problem (\ref{eq22}) has at least one solution for any $\lambda\in (0, \lambda_f)$.}
\vskip 0.1in
{\bf Proof:}\ Using the notation given in Section 1, we have $g(x, \omega)=\omega^p+\lambda f(x)$. To complete the proof of Theorem 2.1, we only need to show that $g(x, \omega)$ verifies {\bf (G1)} and {\bf (G2)}. To verify {\bf (G1)} , we let $\phi(x)$ be the solution of problem (\ref{eq21}). Taking $\varphi=\lambda\phi$, then we have
\begin{equation}\label{eq23}
\left\{
\begin{array}{ll}
-\Delta\varphi=\lambda f(x)\leq\varphi^p+\lambda f(x) & x\in\Omega,\\
\varphi\geq 0 & x\in\Omega,\\
\varphi=0 & x\in\partial\Omega.
\end{array}
\right.
\end{equation}

Let $u(x)$ be the solution of the following problem
\begin{equation}\label{42}
\left\{
\begin{array}{ll}
-\Delta u=1 & x\in\Omega,\\
u=0 & x\in\partial\Omega.
\end{array}
\right.
\end{equation}
Choosing $M_0>0$ so small that
$$M_0>M_0^p\max\limits_{x\in\Omega} u^p(x)+M_0^p\max\limits_{x\in\Omega}|f(x)|,$$
and setting $\psi(x)=M_0u(x)$, we can easily check that
\begin{equation}\label{43}
\left\{
\begin{array}{ll}
-\Delta\psi=M_0\geq\psi^p+\lambda f(x) & x\in\Omega\\
\psi=0 & x\in\partial\Omega
\end{array}
\right.
\end{equation}
for any $\lambda\in (0, M_0^p)$.

Taking (\ref{eq23}) and (\ref{43}) into account, we infer from the strong comparison principle for Laplace operator that
\begin{equation}\label{44}
\varphi(x)<\psi(x)\ \ \mbox{for}\ \ x\in\Omega\ \ \mbox{and}\ \ \lambda\in (0, M_0^p).
\end{equation} 
Let $\lambda_f=M_0^p$. Then for any $\lambda\in(0, \lambda_f)$, we have $0\leq\varphi<\psi$, and
\begin{equation}\label{eq24}
\left\{\begin{array}{ll}
-\Delta\varphi\leq\varphi^p+\lambda f(x) & x\in\Omega,\\
-\Delta\psi\geq \psi^p+\lambda f(x) & x\in\Omega,\\
\varphi=\psi=0 & x\in\partial\Omega.
\end{array}
\right.
\end{equation}
Therefore, {\bf (G1)} is satisfied. Noting that for any $0<\beta\leq 1$ and $\omega\in [\beta\varphi, \psi]$ there hold
$$-\Delta\varphi=\lambda f(x)\leq\omega^p+\lambda f(x)\leq\psi^p+\lambda f(x)\leq-\Delta\psi,$$
we can see that {\bf (G2)} is satisfied with $\alpha=0$. Consequently, the conclusion of Theorem 2.1 follows from Theorem 1.1 and the strong comparison principle of Laplace operator.
\vskip 0.1in
The second example is a Dirichlet problem of Kirchhoff-Carrier type equation involving cocave-convex nonlinearity.
\vskip 0.1in
{\bf Example 2:}\ Assume $c, d>0$ be constants, and $0<q<1<p$. For parameter $\mu>0$, we consider the problem
\begin{equation}\label{eq25}
\left\{
\begin{array}{ll}
-(1+c\|u\|_2^2+d\|\nabla u\|_2^2)\Delta u=\mu u^q+u^p & x\in\Omega,\\
u>0 & x\in\Omega,\\
u=0 & x\in\partial\Omega.
\end{array}
\right.
\end{equation}
By making use of Theorem 1.1, we can prove the following Theorem.
\vskip 0.1in
{\em {\bf Theorem 2.2}\ There exists a positive number $\mu_0$ such that problem (\ref{eq25}) has at least one solution for any $\mu\in (0, \mu_0)$.}
\vskip 0.1in
{\bf Proof:}\ In this example, we have $g(x,\omega)=\mu\omega^q+\omega^p$. To verify {\bf (G1)}, we let $\lambda_1(\Omega)$ be the first eigenvalue of the following eigenvalue problem
\begin{equation}\label{eq26}
\left\{
\begin{array}{ll}
-\Delta\phi=\lambda\phi & x\in\Omega,\\
\phi=0 & x\in\partial\Omega,
\end{array}
\right.
\end{equation}
and denote by $\phi_1(x)$ the positive first eigenfunction which is normalized so that $\max\limits_{x\in\Omega}\phi_1(x)=1$. Since $0<q<1$, we can choose $M_0$ so small that
$$\lambda_1(\Omega)M_0\leq\mu M_0^q.$$
Set $\varphi(x)=M_0\phi_1(x)$. Then, we have
\begin{equation}\label{eq27}
\left\{
\begin{array}{ll}
-\Delta\varphi=M_0\lambda_1(\Omega)\phi_1\leq\mu M_0^q\phi_1^q\leq\mu\varphi^q & x\in\Omega,\\
\varphi=0 & x\in\partial\Omega.
\end{array}
\right.
\end{equation}

Let $\xi(x)$ be the solution of the following problem
\begin{equation}\label{eq28}
\left\{
\begin{array}{ll}
-\Delta\xi=1 & x\in\Omega,\\
\xi=0 & x\in\partial\Omega.
\end{array}
\right.
\end{equation}
Since $p>1$, we can choose $M_1$ so small that
$$M_1>M_1^p\max\limits_{x\in\Omega}\xi^p(x).$$
Setting $\psi(x)=M_1\xi(x)$, and $\mu_0=\frac{M_1-M_1^p\max\limits_{x\in\Omega}\xi^p(x)}{M_1^q\max\limits_{x\in\Omega}\xi^q(x)}$, then, for any $\mu\in (0, \mu_0)$, we have
$$-\Delta\psi=M_1>\mu M_1^q\max\limits_{x\in\Omega}\xi^q(x)+M_1^p\max\limits_{x\in\Omega}\xi^p(x)\geq\mu\psi^q+\psi^p.$$
If necessary, we can choose $M_0$ even more small so that $0<\varphi(x)\leq\psi(x)$. Therefore, for any $\mu\in (0, \mu_0)$, we have 
\begin{equation}\label{eq29}
\left\{
\begin{array}{ll}
-\Delta\varphi\leq\mu\varphi^q+\varphi^p & x\in\Omega,\\
-\Delta\psi\geq\mu\psi^q+\psi^p & x\in\Omega,\\
0<\varphi\leq\psi & x\in\Omega,\\
\varphi=\psi=0 & x\in\partial\Omega.
\end{array}
\right.
\end{equation} 
This implies that {\bf (G1)} hold. Noting that for any $0<\beta\leq 1$ and $\omega\in[\beta\varphi, \psi]$, we can infer from (\ref{eq27}) and (\ref{eq29}) that
$$-\beta^q\Delta\varphi\leq\mu\beta^q\varphi^q\leq\mu\omega^q+\omega^p\leq\mu\psi^q+\psi^p\leq-\Delta\psi.$$
This implies that {\bf (G2)} is valid for $\alpha=q$. Moreover, we can easily see that
$$\lim\limits_{n\to +\infty}\sum\limits_{k=0}^{n-1}q^k=\frac{1}{1-q}.$$
Therefore, we can infer from Theorem 1.1 that problem (\ref{eq25}) has at least one solution $u(x)$ with property 
$$0<(\frac{1}{1+c\|\psi\|_2^2+d\|\nabla\psi\|_2^2})^{\frac{1}{1-q}}\varphi(x)\leq u(x)\leq\psi(x)\ \ \mbox{for}\ \ x\in\Omega.$$
This completes the proof of Theorem 2.2.
\vskip 0.1in
In the above examples, the operator or nonlinear term have more or less monotonicity property. Here, we give an example whose operator and nonlinear term are all not monotone.
\vskip 0.1in
{\bf Example 3:}\ Assume that $f(x)$ satisfies the condition given in Example 1. We consider the following problem
\begin{equation}\label{eq210}
\left\{
\begin{array}{ll}
-(1+d\sin^2(\|\nabla u\|_2))\Delta u=\sin^2(u)+f(x) & x\in\Omega,\\
u>0 & x\in\Omega,\\
u=0 & x\in\partial\Omega.
\end{array}
\right.
\end{equation}
By making use of Theorem 1.1, we prove the following result.
\vskip 0.1in
{\em {\bf Theorem 2.3}\ For any $f(x)$ satisfying the condition given in Example 1, Problem (\ref{eq210}) has at least one solution.}
\vskip 0.1in
{\bf Proof:}\ Let $\varphi(x)$ be the solution of problem (\ref{eq21}). Then we have
$$-\Delta\varphi=f(x)\leq\sin^2(\varphi)+f(x)\ \ \mbox{for}\ \ x\in\Omega.$$
Let $\xi(x)$ be the solution of the following problem
\begin{equation}\label{eq211}
\left\{
\begin{array}{ll}
-\Delta\xi=1 & x\in\Omega,\\
\xi=0 & x\in\partial\Omega.
\end{array}
\right.
\end{equation}
Choosing $M_0>1+\max\limits_{x\in\Omega}|f(x)|$, then $\psi(x)=M_0\xi(x)$ satisfies
$$-\Delta\psi=M_0>\sin^2(\psi)+f(x).$$
Moreover, by strong comparison principle we still have
$$0\leq\varphi(x)<\psi(x)\ \ \mbox{for}\ \ x\in\Omega.$$
Therefore, {\bf (G1)} is satisfied with the above determined $\varphi$ and $\psi$. Noting that for any $0<\beta\leq 1$ and $\omega\in [\beta\varphi, \psi]$ we have
$$-\Delta\varphi=f(x)\leq\sin^2(\omega)+f(x)\leq 1+f(x)<M_0=-\Delta\psi,$$
{\bf (G2)} is satisfied with $\alpha=0$. Consequently, it follows from Theorem 1.1 that problem (\ref{eq210}) has at least one solution $u(x)$ with property $0\leq\varphi(x)\leq u(x)\leq\psi(x)$ for any $x\in\Omega$. Finally, the positivity of $u(x)$ follows from the strong comparison principle of Laplace operator. This completes the proof of Theorem 2.3.
\vskip 0.1in
Almost all $g(x, u)$ in the above examples have definite sign with respect to variable $u$. The following problem provides an example which permits $g(x, u)$ changing sign with respect to variable $u$.
\vskip 0.1in
{\bf Example 4}\ Assume that $0<q<1$ and $d>0$. we consider the following problem
\begin{equation}\label{eq212}
\left\{\begin{array}{ll}
-(1+d\|\nabla u\|_2^2)\Delta u=\mu u^q+\pi\sin^3(u)& x\in\Omega,\\
u>0 & x\in\Omega,\\
u=0 & x\in\partial\Omega.
\end{array}
\right.
\end{equation} 
Using Theorem 1.1, we can conclude the following result
\vskip 0.1in
{\em {\bf Theorem 2.4}\ There exists a positive number $\mu_0$ such that problem (\ref{eq212}) has at least one positive solution for any $\mu\in (0, \mu_0)$.}
\vskip 0.1in
{\bf Proof:}\ Let $\lambda_1(\Omega)$ be the first eigenvalue of the following eigenvalue problem
\begin{equation}\label{eq213}
\left\{
\begin{array}{ll}
-\Delta\phi=\lambda\phi & x\in\Omega,\\
\phi=0 & x\in\partial\Omega,
\end{array}
\right.
\end{equation}
and denote by $\phi_1(x)$ the positive first eigenfunction which is normalized so that $\max\limits_{x\in\Omega}\phi_1(x)=1$. Since $0<q<1$, we can choose $M_0$ so small that
$$\lambda_1(\Omega)M_0\leq\mu M_0^q\ \ \mbox{and}\ \ \sin(M_0\phi_1)\geq 0.$$
Set $\varphi(x)=M_0\phi_1(x)$. Then, we have
\begin{equation}\label{eq214}
\left\{
\begin{array}{ll}
-\Delta\varphi=M_0\lambda_1(\Omega)\phi_1\leq\mu M_0^q\phi_1^q\leq\mu\varphi^q & x\in\Omega,\\
\varphi=0 & x\in\partial\Omega.
\end{array}
\right.
\end{equation}

Let $\xi(x)$ be the solution of the following problem
\begin{equation}\label{eq215}
\left\{
\begin{array}{ll}
-\Delta\xi=1 & x\in\Omega,\\
\xi=0 & x\in\partial\Omega.
\end{array}
\right.
\end{equation}
Choosing $M_1$ so small that
$$M_1>\pi M_1^3\max\limits_{x\in\Omega}\xi^3(x)\ \ \mbox{and}\ \ 0<M_1\xi(x)\leq\frac{\pi}{2}.$$
Setting $\psi(x)=M_1\xi(x)$, and $\mu_0=\frac{M_1-M_1^3\max\limits_{x\in\Omega}\xi^3(x)}{M_1^q\max\limits_{x\in\Omega}\xi^q(x)}$, then, for any $\mu\in (0, \mu_0)$, we have
$$-\Delta\psi=M_1>\mu M_1^q\max\limits_{x\in\Omega}\xi^q(x)+\pi M_1^3\max\limits_{x\in\Omega}\xi^3(x)\geq\mu\psi^q+\pi\sin^3(\psi).$$
If necessary, we can choose $M_0$ even more small so that $0<\varphi(x)\leq\psi(x)$. Therefore, for any $\mu\in (0, \mu_0)$, we have 
\begin{equation}\label{eq216}
\left\{
\begin{array}{ll}
-\Delta\varphi\leq\mu\varphi^q+\pi\sin^3(\varphi) & x\in\Omega,\\
-\Delta\psi\geq\mu\psi^q+\pi\sin^3(\psi) & x\in\Omega,\\
0<\varphi\leq\psi & x\in\Omega,\\
\varphi=\psi=0 & x\in\partial\Omega.
\end{array}
\right.
\end{equation} 
This implies that {\bf (G1)} hold. Noting that $0<\varphi(x)\leq\psi\leq\frac{\pi}{2}$, we can infer from (\ref{eq214}) and (\ref{eq216}) that
$$-\beta^q\Delta\varphi\leq\mu\beta^q\varphi^q\leq\mu\omega^q+\pi\sin^3(\omega)\leq\mu\psi^q+\pi\sin^3(\psi)\leq-\Delta\psi$$
for any $0<\beta\leq 1$ and $\omega\in[\beta\varphi, \psi]$. This implies {\bf (G2)} is valid for $\alpha=q$. Moreover, we can easily see that
$$\lim\limits_{n\to +\infty}\sum\limits_{k=0}^{n-1}q^k=\frac{1}{1-q}.$$
Therefore, we can infer from Theorem 1.1 that problem (\ref{eq212}) has at least one solution $u(x)$ with property 
$$0<(\frac{1}{1+d\|\nabla\psi\|_2^2})^{\frac{1}{1-q}}\varphi(x)\leq u(x)\leq\psi(x)\leq\frac{\pi}{2}\ \ \mbox{for}\ \ x\in\Omega.$$
This completes the proof of Theorem 2.4.

\end{document}